\documentclass[10pt]{article}
\usepackage[utf8]{inputenc}
\usepackage[english]{babel}
\usepackage{amsmath,amsthm,amsfonts}
\usepackage{graphicx} 
\usepackage{bm}
\usepackage{tikz}
\usepackage{multirow}
\usepackage{color}

\newtheorem{remark}{Remark}

\graphicspath{ {./} }

\begin{document}

\title{Generalized Multiscale Finite Element method for multicontinua unsaturated flow problems in fractured porous media}

\author{
Denis Spiridonov
\thanks{Multiscale model reduction laboratory, North-Eastern Federal University, Yakutsk, Republic of Sakha (Yakutia), Russia, 677980. 
Email: {\tt d.stalnov@mail.ru}.}
\and
Maria Vasilyeva \thanks{Institute for Scientific Computation, Texas A\&M University, College Station, TX 77843-3368 \& Department of Computational Technologies, North-Eastern Federal University, Yakutsk, Republic of Sakha (Yakutia), Russia, 677980. 
Email: {\tt vasilyevadotmdotv@gmail.com}.}
\and
Eric T. Chung \thanks{Department of Mathematics,
The Chinese University of Hong Kong (CUHK), Hong Kong SAR. 
Email: {\tt tschung@math.cuhk.edu.hk}.}
}

\maketitle

\begin{abstract}
In this paper, we present a multiscale method for simulations of the multicontinua unsaturated flow problems in heterogeneous fractured porous media.
The mathematical model is described by the system of Richards equations for each continuum that coupled by the specific transfer term.
To illustrate the idea of our approach, we consider a dual continua background model with discrete fractures networks that generalized as a multicontinua model for unsaturated fluid flow in the complex heterogeneous porous media. We present fine grid approximation based on the finite element method and Discrete Fracture Model (DFM) approach.
In this model, we construct an unstructured fine grid that take into account a complex fracture geometries for two and three dimensional formulations. Due to construction of the unstructured grid, the fine grid approximation leads to the very large system of equations. For reduction of the discrete system size, we develop a multiscale method for coarse grid approximation of the coupled problem using Generalized Multiscale Finite Element Method (GMsFEM). In this method, we construct a coupled multiscale basis functions that used to construct highly accurate coarse grid approximation.
The multiscale method allowed us to capture detailed interactions between multiple continua.
The adaptive approach is investigated, where we consider two approaches for multiscale basis functions construction: (1) based on the spectral characteristics of the local problems and (2) using simplified multiscale basis functions.
We investigate accuracy of the proposed method for the several test problems in two and three dimensional formulations.
We present a comparison of the relative error for different number of basis functions and for adaptive approach. Numerical results illustrate that the presented method provide accurate solution of the unsaturated multicontinua problem on the coarse grid with huge reduction of the discrete system size.
\end{abstract}

\section{Introduction}
Prediction of fluid flow in unsaturated soils is an important problem in science and engineering, for example, in agriculture and environmental engineering, oil and gas production, groundwater hydrology etc. 
Mathematical models of the unsaturated filtration processes are described by the Richards equation \cite{b1,b2,b4}, where porous media characterized by a complex heterogeneous rock properties, complex fracture distribution, multiple scales and high contrast of the properties. 
Due to the high permeability of the fractures, they have a significant effect on the flow processes and requires a special approach in the construction of a mathematical model and computational algorithms. 
For example, a hierarchical model is used to describe the multiscale fractures \cite{lee2001hierarchical, li2006efficient, chung2017coupling, li2018multiscale}, where flow in the small scale highly connected fracture networks (natural fractures) can be described by the dual continua approach with additional lower dimensional discrete fracture for simulation of the flow in large scale fractures.

In this work, we construct a multicontinua models for unsaturated flow in porous media, that based on the coupled system of Richards equations. We consider fractured porous media, where one of the continua is describes a flow in the fracture networks. For accurate simulations in the fractured and heterogeneous porous media, the unstructured grids should be used to resolve features in the level of mesh construction. Such grids lead to the large number of unknowns of the discrete system and computationally expensive. To solve such problems, a homogenization techniques or multiscale methods are used to reduce size of the discrete systems. 

Several multiscale methods are developed to solve problems with heterogeneous properties, for example, multiscale finite element method (MsFEM) \cite{b5}, multiscale finite volume method (MsFVM) \cite{hajibeygi2008iterative}, heterogeneous multiscale methods (HMM) \cite{weinan2007heterogeneous}, generalized multiscale finite element method (GMsFEM) \cite{CYH2016adaptive}, constraint energy minimizing generalized multiscale finite element method (CEM-GMsFEM)  \cite{chung2018constraint}.  
In \cite{he2009adaptive, he2006multiscale, b3}, the authors present construction of the coarse grid approximation based on the MsFEM  for solving the unsaturated flow problems with heterogeneous coefficients. Upscaling method for the Richards equation  is presented in \cite{chen2005upscaling}.  
Multiscale methods for solution of the flow problems in fractured porous media are presented in \cite{akkutlu2015multiscale, ctene2016algebraic, bosma2017multiscale, ctene2017projection, hkj12}.  
In our previous works, we developed multiscale model reduction techniques based on the Generalized Multiscale Finite Element Method (GMsFEM) for flow in fractured porous media \cite{akkutlu2015multiscale, chung2017coupling, akkutlu2018multiscale}. 
Recently in \cite{chung2018non, vasilyeva2019nonlocal, vasilyeva2019constrained}, we introduced a nonlocal multicontinuum (NLMC) method for problems in fractured  porous media, where we construct multiscale basis functions based on the solution of some local constrained energy minimization problems as in the CEM-GMsFEM.

In this work, we present the Generalized Multiscale Finite Element method \cite{b6,b8} for solution of the unsaturated multicontinua flow problem in fractured heterogeneous porous media. Previously in \cite{chung2017coupling}, we presented construction of the GMsFEM for solution of the saturated flow problem in fractured media, in presented work we extend proposed framework for unsaturated flow case problem in two and three dimensional formulations. 
In GMsFEM, we construct multiscale basis functions that automatically identify each continuum via solution of the local spectral problems  \cite{b10,b11, b12}.
For coupled system of equations in multicontinua models, we solve local coupled system of equations to construct highly accurate basis functions that capture complex interaction between continua in the coarse grid model. 
We study an adaptive approach with simplified basis functions construction and compare with adaptive choosing based on the local eigenvalues of the spectral problem. The simplified approach is effective for the cases, when the fracture networks are known and have simplified geometries. 
To illustrate  robustness and accuracy of the proposed method, we consider a dual continuum background model that coupled with the discrete fracture networks (triple continua model). We present several numerical examples for two and three dimensional test problems.

We organize the paper as follows. In Section 2, we present a mathematical model for unsaturated multicontinua flow problem in fractured porous media as a general multicontinia formulation of the unsaturated flow. In section 3, we consider fine-scale approximation using Discrete Fracture Model. Next in Section 4, we describe the construction of the coarse grid approximation and describe construction of the multiscale basis functions. Furthermore, we present an adaptive approach with simplified multiscale basis functions and using adaptive approach based on the spectral characteristics of the spectral problems. We present numerical results in Sections 5 for two and three dimensional model problems. Finally, we present Conclusion.

\section{Mathematical model}

Mathematical model of the unsaturated flow in porous media described by the Richard's equations. For fractured porous media, we consider a mixed dimensional formulation of the flow problem \cite{martin2005modeling, d2012mixed, formaggia2014reduced, Quarteroni2008coupling}.
Let $\Omega \in \mathcal{R}^d$ is the $d$ - dimensional domain of the porous matrix, where $d = 2,3$. Fracture network is considered as a $(d-1)$ - dimensional (lower dimensional) domain $\gamma \in \mathcal{R}^{d-1}$ due to small thickness of the fractures compared to the domain sizes.
Then for unsaturated flow in fractured porous media, we have following coupled system of equations
\begin{equation}
\label{mm1}
\begin{split}
& \frac{\partial \Theta_m}{\partial t} + \nabla \cdot q_{m} + L_{mf} = f_m, \quad x \in \Omega, \\
& q_{m} = - k_{m}(x, p_m) \nabla (p_{m} + z), \quad x \in \Omega, \\
& \frac{\partial \Theta_f}{\partial t} + \nabla \cdot q_{f} - L_{fm} = f_f, \quad x \in \gamma, \\
& q_{f} = - k_{f}(x, p_f) \nabla (p_{f} + z), \quad x \in \gamma,
\end{split}
\end{equation}
where $q_m$ and $q_f$ are the Darcy velocities in matrix and fractures;
$p_m$ and $p_f$ are the pressure head in matrix and fractures;
$k_m$ and $k_f$ are the unsaturated hydraulic conductivity tensors for matrix and fractures;
$z$ represent the influence of the gravity to the flow processes;
$\Theta_m$ and $\Theta_f$ are the water content for matrix and fracture; and
$f_m$ and $f_f$ refer to source and sink terms.
Here $L_{mf}$ and $L_{fm}$ are the transfer term between matrix-fracture and fracture-matrix, $L_{mf} = \eta_m q_{mf}$, $L_{fm} = \eta_f q_{mf}$, $q_{mf} \approx \sigma_{mf}(x, p_m, p_f) (p_m - p_f)$ and $\int_{\Omega} L_{mf} dx - \int_{\gamma} L_{fm} ds = 0$. 

After substitution the Darcy's Law into the mass conservation equation, we obtain following system of equations for $p_m$ and $p_f$:
\begin{equation}
\label{mm-fm}
\begin{split}
& \frac{\partial \Theta_m}{\partial t}
- \nabla \cdot ( k_{m}(x, p_m) \nabla (p_{m} + z))
+ L_{mf} = f_m, \quad x \in \Omega, \\
& \frac{\partial \Theta_f}{\partial t}
- \nabla \cdot ( k_{f}(x, p_f) \nabla (p_{f} + z))
- L_{fm} = f_f, \quad x \in \gamma.
\end{split}
\end{equation}

Similarly, for the dual continuum background model with lower-dimensional fracture model, we have following system of equations for $p_1$, $p_2$ and $p_f$:
\begin{equation}
\label{mm-tc}
\begin{split}
& \frac{\partial \Theta_1}{\partial t}
- \nabla \cdot ( k_1(x, p_1) \nabla (p_1 + z))
+ L_{12} + L_{1f} = f_1, \quad x \in \Omega, \\
& \frac{\partial \Theta_2}{\partial t}
- \nabla \cdot ( k_2(x, p_2) \nabla (p_2 + z))
- L_{21} + L_{2f} = f_2, \quad x \in \Omega, \\
& \frac{\partial \Theta_f}{\partial t}
- \nabla \cdot ( k_{f}(x, p_f) \nabla (p_{f} + z))
- L_{f1} - L_{f2} = f_f, \quad x \in \gamma.
\end{split}
\end{equation}
where subindices $1, 2$ are related for first and second continua, and subindex $f$ is related to the lower-dimensional fracture model.
System of equations is coupled by the mass exchange terms between continuum and $\int_{\Omega} L_{\alpha f} dx - \int_{\gamma} L_{f\alpha} ds = 0$, $\int_{\Omega} L_{\alpha \beta} dx - \int_{\Omega} L_{\beta \alpha} dx = 0$ for $\alpha, \beta = 1,2$ and suppose following linear relation $L_{12} = L_{21} \approx \sigma_{12}(x, p_1, p_2) (p_1 - p_2)$.

We consider system of equations \eqref{mm-tc} unsaturated flow model with following initial conditions
\begin{equation}
\label{mm-ic}
p_1 = p_2 = p_f = p^0, 
\end{equation}
and boundary conditions
\begin{equation}
\label{mm-bc}
\begin{split}
&p_1 = p_2 = p_f = g, \quad x \in \Gamma_D, \\
&q_1 \cdot n = q_2 \cdot n = q_f \cdot n = 0, \quad x \in \Gamma_N,
\end{split}
\end{equation}
where $\partial \Omega = \Gamma_D \cup \Gamma_N$, $\Gamma_D$ is the top boundary of the computational domain and $\Gamma _N$ denotes left, right and bottom boundaries of the $\Omega$. On the fracture boundaries that inside domain $\Omega$, we set a zero Neumann boundary conditions. 

In general, we have following multicontinuum model:
\begin{equation}
\label{mm-mc}
\begin{split}
& \frac{\partial \Theta_{\alpha}}{\partial t}
- \nabla \cdot ( k_{\alpha}(x, p_{\alpha}) \nabla (p_{\alpha} + z))
+ \sum_{\beta \neq \alpha} L_{\alpha \beta} = f_{\alpha},
\end{split}
\end{equation}
where $\alpha = 1,...,L$ and $L$ is the number of continuum.

As constitutive relations, we use Haverkamp model \cite{b4}
\begin{equation}
\label{ThetaK}
\Theta_{\alpha}(p_{\alpha}) =
\frac{A_{\alpha} (\Theta_{\alpha,s} - \Theta_{\alpha,r})}{A_{\alpha} + |p_{\alpha}|^{B_{\alpha}} } + \Theta_{\alpha,r},
\quad
k_{\alpha}(x, p_{\alpha}) = k_{\alpha, s}(x) \frac{C_{\alpha}}{C_{\alpha} + |p_{\alpha}|^{D_{\alpha}} },
\end{equation}
where $A_{\alpha}$ $B_{\alpha}$, $C_{\alpha}$, $D_{\alpha}$, $\Theta_{\alpha,r}$ and $\Theta_{\alpha,s}$ are the Haverkamp model coefficients, $k_{\alpha,s}(x)$ is the heterogeneous saturated hydraulic conductivity for $\alpha$ continuum.

\section{Fine grid approximation}

For numerical solution, we construct unstructured fine grid that explicitly resolve fractures in the level of mesh. We construct a discrete system based on the finite element method and Discrete Fracture Model (DFM) for fracture networks.  
Let $\mathcal{T}_h$ denote a finite element partition of the domain $\Omega \in \mathcal{R}^d$ and  $\mathcal{T}_h = \bigcup_i K_i$, where $K_i$ is the triangular element for $d = 2$ and tetrahedron for $d = 3$. 
Let $\mathcal{E}_h$ the set of all the faces between the elements $\mathcal{T}_h$ and  $\mathcal{E}_{\gamma} \subset \mathcal{E}_h$ be the subset of all faces that represent fractures. Furthermore, fracture facets $\mathcal{E}_{\gamma}$ represent the lower dimensional fracture grid. 
Then finite element approximation for system of equations \eqref{mm-tc} with boundary conditions \eqref{mm-bc} can be written as follows: find $p = (p_1, p_2, p_f) \in V_1 \times V_2 \times V_f$ such that
\begin{equation} 
\label{fapp-tc}
\begin{split}
m_1 \left(\frac{\partial \Theta_1}{\partial t}, v_1\right) + a_1 (p_1, v_1) 
+ q_{12}(p_1-p_2, v_1) + q_{1f}(p_1-p_f, v_1)  = l(v_1), \\
m_2 \left(\frac{\partial \Theta_2}{\partial t}, v_2\right) + a_2 (p_2, v_2) 
- q_{21}(p_1-p_2, v_2) + q_{2f}(p_2-p_f, v_2)  = l(v_2), \\
m_f \left(\frac{\partial \Theta_f}{\partial t}, v_f\right) + a_f (p_f, v_f) 
- q_{1f}(p_1-p_f, v_f) - q_{2f}(p_2-p_f, v_f)  = l(v_f),  
\end{split}
\end{equation}
for $\forall (v_1, v_2, v_f) \in \hat{V}_1 \times \hat{V}_2 \times \hat{V}_f$ and 
\[
V_1 = V_2 = \{v \in H^1(\Omega): v = g , \, x \in \Gamma_D\}, \quad 
V_f = \{v \in H^1(\gamma): v = g , \, x \in \Gamma_D\},
\]\[
\hat{V}_1 = V_2 = \{v \in H^1(\Omega): v = 0 , \, x \in \Gamma_D\}, \quad 
\hat{V}_f = \{v \in H^1(\gamma): v = 0 , \, x \in \Gamma_D\}.
\]
For bilinear and linear forms, we have
\[
m_{\alpha} \left(\frac{\partial \Theta_{\alpha}}{\partial t}, v_{\alpha}\right) = 
\int_{\Omega}  \frac{\partial \Theta_{\alpha}}{\partial t} \, v_{\alpha} \, dx, \quad 
m_f \left(\frac{\partial \Theta_f}{\partial t}, v_f\right) = 
\int_{\gamma}  \frac{\partial \Theta_f}{\partial t} \, v_f \, ds, 
\]\[
a_{\alpha} (p_{\alpha}, v_{\alpha}) = 
\int_{\Omega} k_{\alpha}(x, p_{\alpha})  \nabla p_{\alpha} \cdot \nabla v_{\alpha} \, dx, \quad 
a_f (p_f, v_f)  =
 \int_{\gamma} k_f(x, p_f)  \nabla p_f \cdot \nabla v_f \, ds, 
\]\[
l(v_{\alpha}) =  \int_{\Omega} f_{\alpha} \, v_{\alpha} \, dx 
+ \int_{\Omega} \frac{\partial k_{\alpha}(x, p_{\alpha})}{\partial z}v_{\alpha} \, dx, \quad
l(v_f) = \int_{\gamma} f_f \, v_f \, ds 
+ \int_{\gamma} \frac{\partial k_f (x, p_f)}{\partial z}v_f \, ds,
\]\[
q_{12}(p_{1}-p_{2}, v_{1}) = q_{21}(p_{1}-p_{2}, v_{2}) = 
\int_{\Omega} \sigma_{12}(p_1 - p_2) \, v_2 \, dx, 
\]\[ 
q_{\alpha f}(p_{\alpha}-p_f, v_{\alpha})  = q_{\alpha f}(p_{\alpha}-p_f, v_f) = 
\int_{\gamma} \sigma_{\alpha f}(\mathcal{P}^{\gamma}(p_{\alpha}) - p_f)  \, v_f \, ds 
\]
for $\alpha, \beta = 1,2$.

Similarly, we can generalize model for multicontinuum case: find $p = (p_1, ..., p_L) \in V_1 \times ... \times V_L$ such that
\begin{equation}
\label{fapp-mc}
\begin{split}
m_{\alpha} \left(\frac{\partial \Theta_{\alpha}}{\partial t}, v_{\alpha}\right)
+ a_{\alpha} (p_{\alpha}, v_{\alpha})
+ \sum_{\beta \neq \alpha}
q_{\alpha \beta}(p_{\alpha}-p_{\beta}, v_{\alpha})
= l(v_{\alpha}), \quad \alpha = 1,...,L,
\end{split}
\end{equation}
for $\forall (v_1, ..., v_L) \in \hat{V}_1 \times ... \times \hat{V}_L$.

To solve the coupled nonlinear system of equations \eqref{fapp-mc}, we apply implicit approximation by time with a Picard iteration scheme.
For $\Theta_{\alpha}^{n+1, m+1}$ we use Taylor series with respect to $p_{\alpha}$
\begin{equation}
\label{fapp-theta}
\begin{split}
\Theta_{\alpha}^{n+1, m+1}
&\approx \Theta_{\alpha}^{n+1, m}
+ \left. \frac{ d \Theta_{\alpha} }{dp_{\alpha}} \right|^{n+1, m} (p_{\alpha}^{n+1, m+1} - p_{\alpha}^{n+1, m}) \\
&= \Theta_{\alpha}^{n+1, m} + C_{\alpha}^{n+1, m} (p_{\alpha}^{n+1, m+1} - p_{\alpha}^{n+1, m}),
\end{split}
\end{equation}
where $C_{\alpha} =  d \Theta_{\alpha}/ dp_{\alpha}$,  superscript $n$ indicate time iteration and $m$ is the nonlinear iteration number.

Therefore, we obtain following variational formulation
\begin{equation}
\label{fapp-mc-nt}
\begin{split}
c^{n+1, m}_{\alpha} &\left( \frac{ p_{\alpha}^{n+1, m+1} - p_{\alpha}^{n}}{\tau}, v_{\alpha} \right) 
+ a^{n+1, m}_{\alpha} (  \nabla p_{\alpha}^{n+1, m+1} \cdot \nabla v_{\alpha} )  \\
&+ \sum_{\beta \neq \alpha} 
q^{n+1, m}_{\alpha \beta}(p^{n+1, m+1}_{\alpha} - p^{n+1, m+1}_{\beta}, v_{\alpha}) 
= l^{n+1, m}(v_{\alpha})
- m_{\alpha} \left( \frac{ \Theta_{\alpha}^{n+1, m} - \Theta_{\alpha}^{n}}{\tau}, v_{\alpha} \right), 
\end{split}
\end{equation}
where $\alpha = 1,...,L$, $\tau$ is the given time step and 
\[
c^{n+1, m}_{\alpha} ( p_{\alpha}, v_{\alpha} ) 
=  \int_{\Omega_{\alpha}}  C^{n+1, m}  p_{\alpha} \, v_{\alpha} \, dx, \quad 
\quad
l^{n+1, m}(v_{\alpha}) 
=  \int_{\Omega_{\alpha}} f_{\alpha} \, v_{\alpha} \, dx 
+ \int_{\Omega_{\alpha}} \frac{\partial k_{\alpha}(x, p^{n+1, m}_{\alpha})}{\partial z}v_{\alpha} \, dx,
\]\[
a^{n+1, m}_{\alpha} (p_{\alpha}, v_{\alpha}) 
= \int_{\Omega_{\alpha}} k_{\alpha}(x, p^{n+1, m}_{\alpha})  \nabla p_{\alpha} \cdot \nabla v_{\alpha} \, dx,
\quad
q^{n+1, m}_{\alpha \beta}(p_{\alpha}-p_{\beta}, v_{\alpha}) 
= \int_{\Omega_{\alpha}} \sigma^{n+1, m}_{\alpha \beta}(p_{\alpha} - p_{\beta}) \, v_{\alpha} \, dx.
\]

For triple continuum model \eqref{fapp-tc} with the Galerkin finite element method, we write the solution as
\[
p_1 = \sum_{i = 1}^{N_f^{\Omega}} p_{1,i} \phi_i, \quad
p_2 = \sum_{i = 1}^{N_f^{\Omega}} p_{2,i} \phi_i, \quad
p_f = \sum_{i = 1}^{N_f^{\gamma}} p_{f,i} \psi_i,
\]
where $\phi_i$ are the standard $d$ - dimensional linear element basis functions defined on $\mathcal{T}_h$, and $\psi_i$ are $(d-1)$ - dimensional linear basis function on the fracture mesh, $N_f^{\Omega}$ and $N_f^{\gamma}$ denotes the number of the nodes for $\mathcal{T}_h$ and for fracture mesh.

We write approximation for \eqref{fapp-tc}, in the matrix form as follows for $p = (p_1, p_2, p_f)^T$
\begin{equation}
\label{fapp-tc-m}
C^{n+1, m} \frac{p^{n+1, m+1} - p^{n}}{\tau} + A^{n+1, m} p^{n+1, m+1} = F^{n+1, m},
\end{equation}
where
\[
C^{n+1, m} = 
\begin{pmatrix}
C^{n+1, m}_1 & 0 & 0 \\
0 & C^{n+1, m}_2 & 0 \\
0 & 0 & C^{n+1, m}_f
\end{pmatrix}, \quad 
F^{n+1, m} = 
\begin{pmatrix}
F^{n+1, m}_1 \\ 
F^{n+1, m}_2 \\ 
F^{n+1, m}_f
\end{pmatrix},
\]\[
A^{n+1, m} = 
 \begin{pmatrix}
A^{n+1, m}_1+Q^{n+1, m}_{12}+Q^{n+1, m}_{1f} & -Q^{n+1, m}_{12} & -Q^{n+1, m}_{1f} \\
 -Q^{n+1, m}_{12} & A^{n+1, m}_2+Q^{n+1, m}_{12}+Q^{n+1, m}_{2f} &  -Q^{n+1, m}_{2f} \\
 -Q^{n+1, m}_{1f} &  -Q^{n+1, m}_{2f} & A^{n+1, m}_f + Q^{n+1, m}_{1f} +  Q^{n+1, m}_{2f}
\end{pmatrix}. 
\]
This fine-scale discretization yields matrices of the size $N_f \times N_f$ with  $N_f = N_f^{\Omega} + N_f^{\Omega} + N_f^{\gamma}$.

\begin{remark}
In this paper, for simplification of the matrix construction, we use a modified DFM approach and consider the case when $\sigma_{2f} = 0$. We assume that 
$p_1= p_f$  and using superposition principle \cite{efendiev2015hierarchical, akkutlu2015multiscale, chung2017coupling}, we eliminate $p_f$ from equations \eqref{fapp-tc-m} and obtain following coupled system of equations for $p = (p_1, p_2)^T$
\[
C^{n+1, m} \frac{p^{n+1, m+1} - p^{n}}{\tau} + A^{n+1, m} p^{n+1, m+1} = F^{n+1, m},
\]
where
\[
C^{n+1, m} = 
\begin{pmatrix}
C^{n+1, m}_1+C^{n+1, m}_f & 0 \\
0 & C^{n+1, m}_2 
\end{pmatrix}, \quad 
F^{n+1, m} = 
\begin{pmatrix}
F^{n+1, m}_1 + F^{n+1, m}_f \\ 
F^{n+1, m}_2 
\end{pmatrix},
\]\[
A^{n+1, m} = 
 \begin{pmatrix}
A^{n+1, m}_1 + A^{n+1, m}_f+Q^{n+1, m}_{12} & -Q^{n+1, m}_{12} \\
 -Q^{n+1, m}_{12} & A^{n+1, m}_2+Q^{n+1, m}_{12}
\end{pmatrix}. 
\]
with matrices of the size $N_f \times N_f$ with  $N_f = N_f^{\Omega} + N_f^{\Omega}$.
\end{remark}

\section{Coarse grid approximation using GMsFEM}

Let $\mathcal{T}_H$ is the coarse grid for computational domain $\Omega$ and $\omega_i$ is the local domains, where $i = 1,...,N_c^{\Omega}$ and $N_c^{\Omega}$ is the number of coarse grid nodes. 
For construction of the coarse grid solver, we use a Generalized Multiscale Finite Element Method (GMsFEM) \cite{b12, eh09}.   
For the construction of the multiscale space for coarse grid approximation, we solve spectral problems in each local domain $\omega_i$, $i = 1,...,N_c^{\Omega}$ to identify the most important characteristics of the problem.  Let $p_{C}=\sum_{i, k} p_{k}^i \psi_{k}^i (x)$, where $\psi_{k}^{i}$ are the multiscale basis functions that supported in local domain $\omega_i$, and the index $k$ represents the numbering of the basis functions, $k = 1,..,M_i$. 

\begin{figure}[h!]
\begin{center}
\includegraphics[width=0.99\linewidth]{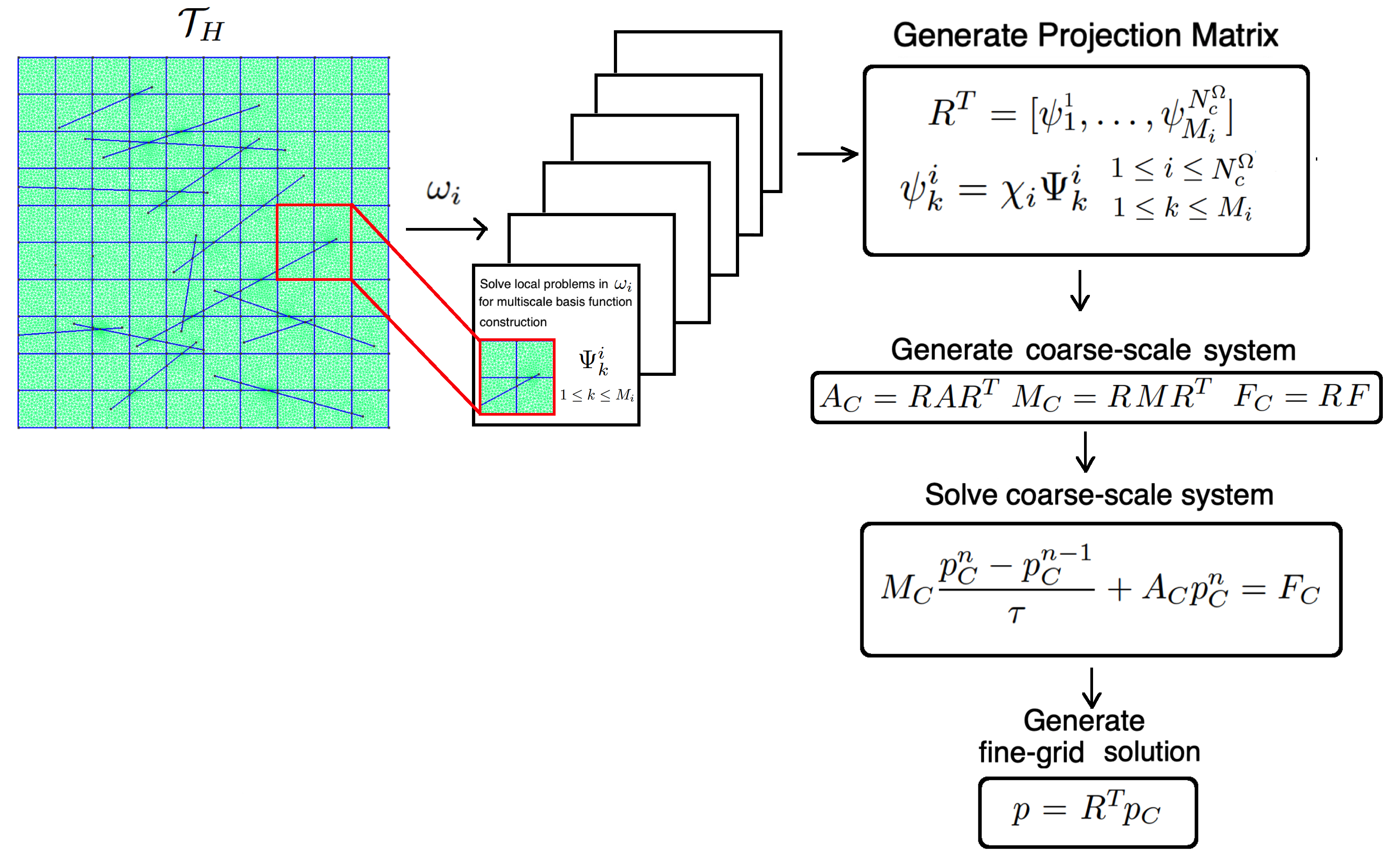} 
\end{center}
\caption{Illustration of the GMsFEM algorithm}
\label{sch-alg}
\end{figure}

In GMsFEM, we have offline and online steps  (see Figure \ref{sch-alg}): 
\begin{itemize}
\item Offline steps: 
\begin{itemize}
\item Coarse grid generation and local domains construction, $\omega_i$.
\item Construction of the local multiscale basis functions, $\psi_{k}^{i}$  in local domain $\omega_i$  ($k = 1,..,M_i$). 
\item Generation of the projection matrix, $R^T = (\psi_1^1,...,\psi_{M_1}^1, ..., \psi_1^{N^{\Omega}_c},...,\psi_{M_{N^{\Omega}_c}}^{N^{\Omega}_c})$.
\end{itemize}
\item Online steps: For each time step and nonlinear iteration:
\begin{itemize}
\item Construct coarse grid system using projection matrix.
\item Solve coarse grid system and reconstruct fine-scale solution. 
\end{itemize}
\end{itemize}
Here, for a given configuration  of  heterogeneous properties and fracture geometry, we precompute a multiscale basis functions on the offline step and generate projection matrix. On the online  step, we project nonlinear matrices into the multiscale space and solve coarse grid systems.  
Note that, we focus method on the structured coarse grids in this work for illustration of the method, but in general, the method in presented form can be extended  to unstructured coarse grids. 

In \eqref{ThetaK}, we consider the case, when the nonlinearity and heterogeneity of $k_{\alpha}(x, p_{\alpha})$ is separable, 
\[
k_{\alpha}(x, p_{\alpha}) = k_{\alpha, r}(p_{\alpha}) k_{\alpha, s}(x).
\] 
Therefore, we can use a linear multiscale space and  precompute multiscale basis functions.  

For general multicontinuum case, let  $p = (p_1, ..., p_L) \in V_1 \times,...,\times V_L$  and
\begin{equation} 
\label{capp-op}
\begin{split}
a^s_{\alpha} (p_{\alpha}, v_{\alpha})  + \sum_{\beta \neq \alpha} 
q^s_{\alpha \beta}(p_{\alpha}-p_{\beta}, v_{\alpha}) 
= 0, \quad \alpha = 1,...,L, 
\end{split}
\end{equation}
where
\[
a^s_{\alpha} (p_{\alpha}, v_{\alpha}) = 
\int_{\Omega_{\alpha}} k_{\alpha, s}(x)  \nabla p_{\alpha} \cdot \nabla v_{\alpha} \, dx, \quad 
q^s_{\alpha \beta}(p_{\alpha}-p_{\beta}, v_{\alpha}) =  
\int_{\Omega_{\alpha}} \sigma^s_{\alpha \beta}(p_{\alpha} - p_{\beta}) \, v_{\alpha} \, dx
\]
and $V_{\alpha} = H^1(\Omega_{\alpha})$.

Next, we describe construction of the coupled mutiscale basis functions for the multicontinua flow problem.  We start with the construction of the snapshot space in the local domain $\omega_i$.  The constructed local snapshots contain the information about local heterogeneities. After, the multiscale spaces are obtained from the  snapshot spaces by a dimension reduction
via local spectral problems \cite{chung2017coupling, akkutlu2015multiscale, akkutlu2017multiscale, efendiev2015hierarchical}.

The snapshot space  is constructed by the solution of the local problems in $\omega_i$ with all possible boundary conditions. 
Let 
$V_{\alpha} = \{ v \in H^1(\omega_i): v = \delta_l  \, \text{on} \,  \partial \omega_i  \}$ ($\alpha = 1,...,L$), 
we find $\phi^{i,l} = (\phi^{i,l}_1, ..., \phi^{i,l}_L) \in V_1 \times ... \times V_L$ such that
\begin{equation}
\label{capp-snap}
a^s_{\alpha} (\phi^{i,l}_{\alpha}, v_{\alpha})  + \sum_{\beta \neq \alpha} 
q^s_{\alpha \beta}(\phi^{i,l}_{\alpha} - \phi^{i,l}_{\beta}, v_{\alpha}) 
= 0, \quad \alpha = 1,...,L, 
\end{equation}
for $\forall v_{\alpha} \in V^0_{\alpha} = \{ v \in H^1(\omega_i): v = 0 \, \text{on} \, \partial \omega_i \}$. 
Here  $\delta_j$ is the function, which takes the value one at node $x_l \in \partial \omega_i$ and zero elsewhere, $l = 1,...,J^i$ ($J^i$ is the number of nodes on the boundary of $\omega_i$).  
Therefore, we define  
\[
V^{snap,i} = \text{span}\{ \phi^{i,l}_1, ...,\phi^{i,l}_{J_i}\}, \quad \text{ and } \quad 
R^T_{snap,i} = ( \phi^{i,l}_1, ...,\phi^{i,l}_{J_i} ). 
\]

For multiscale basis functions $\Psi^i =  (\Psi^i_1, ..., \Psi^i_L)^T$ in $\omega_i$, we solve the local spectral problem on the snapshot space for 
\begin{equation} 
\label{capp-spec}
A_{snap} \Psi^{snap,i} = \lambda S_{snap} \Psi^{snap,i}, 
\end{equation}  
where 
$\Psi^i_{\alpha,k} = R_{snap,i}  \Psi^{snap,i}_{\alpha,k}$. 
We choose the smallest $M_{i}$ eigenvalues and use them for the construction of multiscale basis functions for $k = 1,2..,M_i$. 

Here
\[
A^{snap}  = R_{snap,i} A^s R_{snap,i}^T, \quad 
A^{snap}  = R_{snap,i} S^s R_{snap,i}^T, 
\]\[
S^s = 
\begin{pmatrix}
S^s_1 & 0 & ... & 0 \\
0 & S^s_2 & ... & 0 \\
... & ... & ... & ... \\
0 & 0 & ... & S^s_L
\end{pmatrix}, \quad 
A^s = 
 \begin{pmatrix}
A^s_1+\sum Q^s_{1\beta} & -Q^s_{12}  & ... & -Q^s_{1L} \\
 -Q^s_{12}  & A^s_2+\sum Q^s_{2\beta} & ... & -Q^s_{2L} \\
... & ... & ... &  ... \\
 -Q^s_{1L} &  -Q^s_{2L} & ... & A^s_L + \sum Q^s_{L\beta} 
\end{pmatrix}. 
\]
and 
$A^s_{\alpha} = \{ a^s_{\alpha,lk} \}$, 
$S^s_{\alpha} = \{ s^s_{\alpha,lk} \}$ for $\alpha = 1,...,L$ with
\[
a^s_{\alpha,lk} = \int_{\omega_i} k_{\alpha,s}(x) \nabla \varphi_{\alpha,l} \cdot \nabla \varphi_{\alpha,k} dx, \quad 
s^s_{\alpha,lk}  = \int_{\omega_i} k_{\alpha,s}(x) \, \varphi_{\alpha,l} \, \varphi_{\alpha,k} \, dx.
\]
where $\varphi_{\alpha,l}$ are the linear basis functions.
Note that, the adaptive approach can be applied based on the calculated eigenvalues $\lambda_k$ of the spectral problem \eqref{capp-spec}.


\begin{remark}[Simplified multiscale basis functions]
For the simplified heterogeneities cases, where fracture networks are known and have simplified geometries, the simplified approach for multiscale basis functions construction can be implemented. Simplified multiscale basis functions are constructed without solution of the local spectral problems on the snapshot spaces and defined by a direct setting specific boundary conditions in local problem (see \cite{chung2017coupling} for detailed explanations).
\end{remark}

For construction of the continuous multiscale space, let $\chi_i$ are the linear basis functions in the local domains $\omega_i$ on the coarse grid.
The multiscale space is defined as
\[
V_c = span(
\psi^1_1, ... ,\psi^1_{M_1}, ....,
\psi^{N^{\Omega}_c}_1, ...,\psi^{N^{\Omega}_c}_{M_{N^{\Omega}_c}} ),
\]
with following projection matrix
\[
R^T = ( 
\psi^1_1, ... ,\psi^1_{M_1}, ...., 
\psi^{N^{\Omega}_c}_1, ...,\psi^{N^{\Omega}_c}_{M_{N^{\Omega}_c}} ). 
\] 
where $\psi^i_k = \chi_i \Psi^i_k$ and $\chi_i$ is the standard linear partition of unity functions.

Therefore, we have following coarse grid approximation for the multicontinua problem
\begin{equation}
\label{eq:mfa}
C_c^{n+1, m} \frac{ p_c^{n+1, m+1} -  p_c^{n+1, m} }{\tau} + 
A_c^{n+1, m} p_c^{n+1, m+1} = F_c^{n+1, m},
\end{equation}
where 
$C_c^{n+1, m} = R C^{n+1, m}  R^T$, 
$A_c^{n+1, m} = R A^{n+1, m}  R^T$,  
$ F_c^{n+1, m} = R  F^{n+1, m}$ and the fine grid solution can be reconstructed, $p_{ms} = R^T p_c$.

\section{Numerical results}

\begin{figure}
\centering
\includegraphics[width=0.9\linewidth]{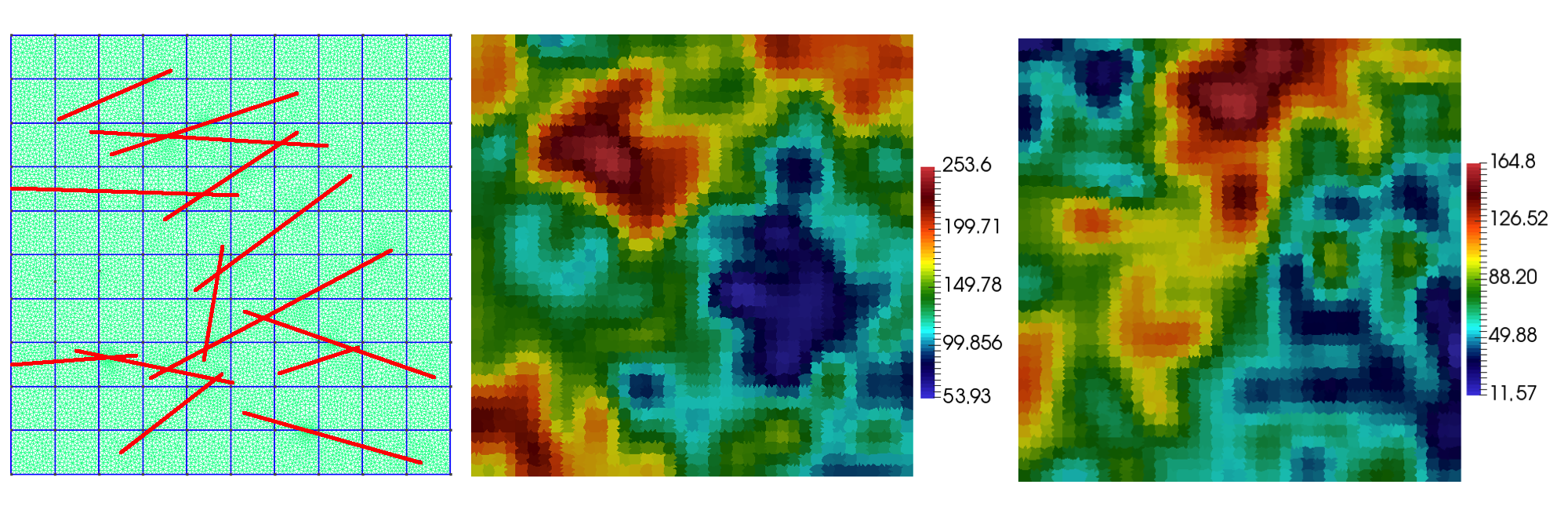}
\caption{Computational grids and heterogeneous properties for \textit{Test 2} and or \textit{Test 3} (two-dimensional problem). 
Left: coarse grid (blue color),  fine grid (green color) and fracture geometry (red color). 
Middle: Heterogeneous coefficient $k_{1,s}(x)$ for \textit{Test 3}. 
Right:  Heterogeneous coefficient $k_{2,s}(x)$ for \textit{Test 3}  ($k_{m,s}(x)$ for \textit{Test 2}). }
\label{geom2d}
\end{figure}

\begin{figure}[h!]
\centering
\includegraphics[width=0.6\linewidth]{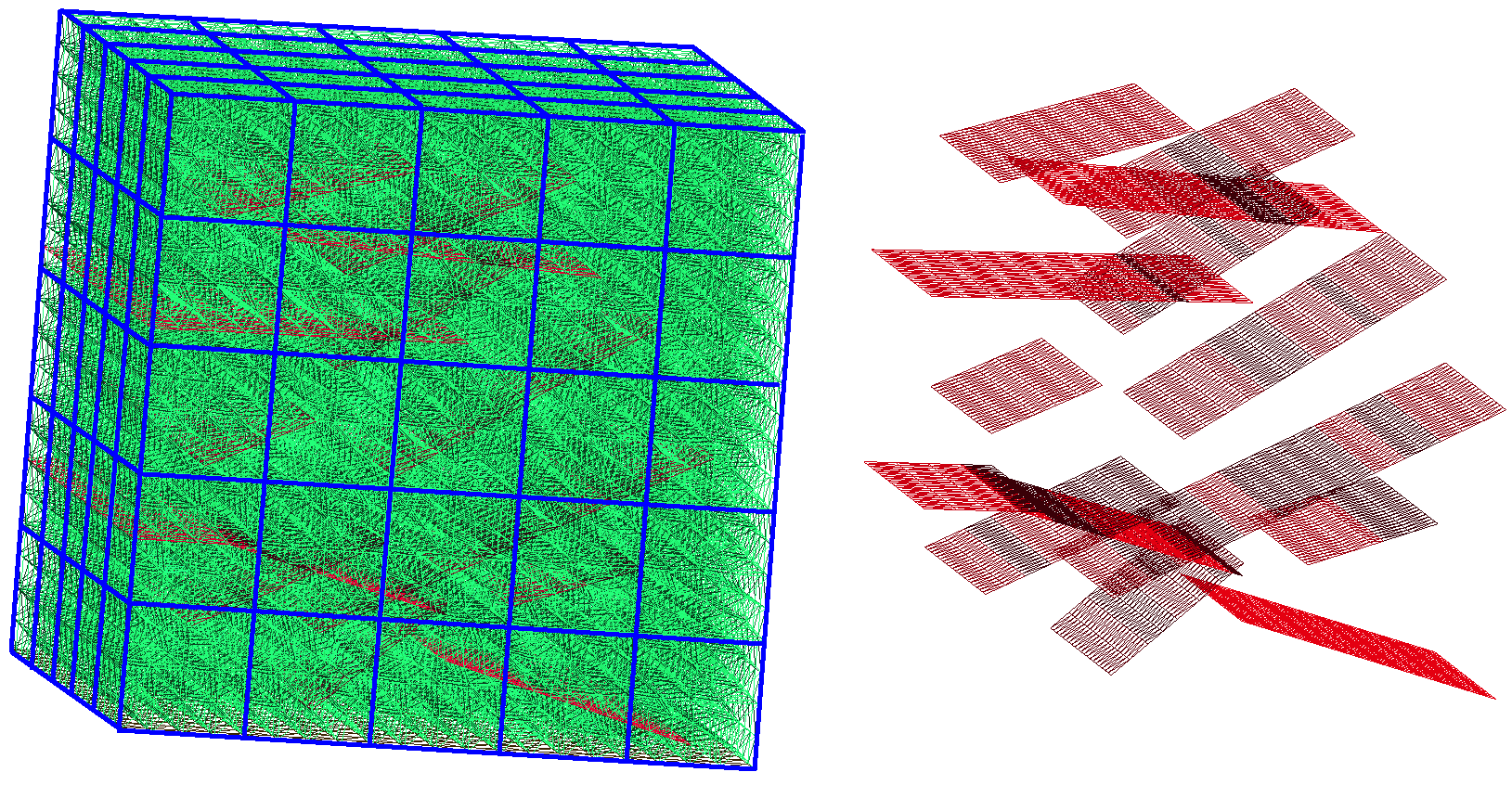} 
\includegraphics[width=0.8\linewidth]{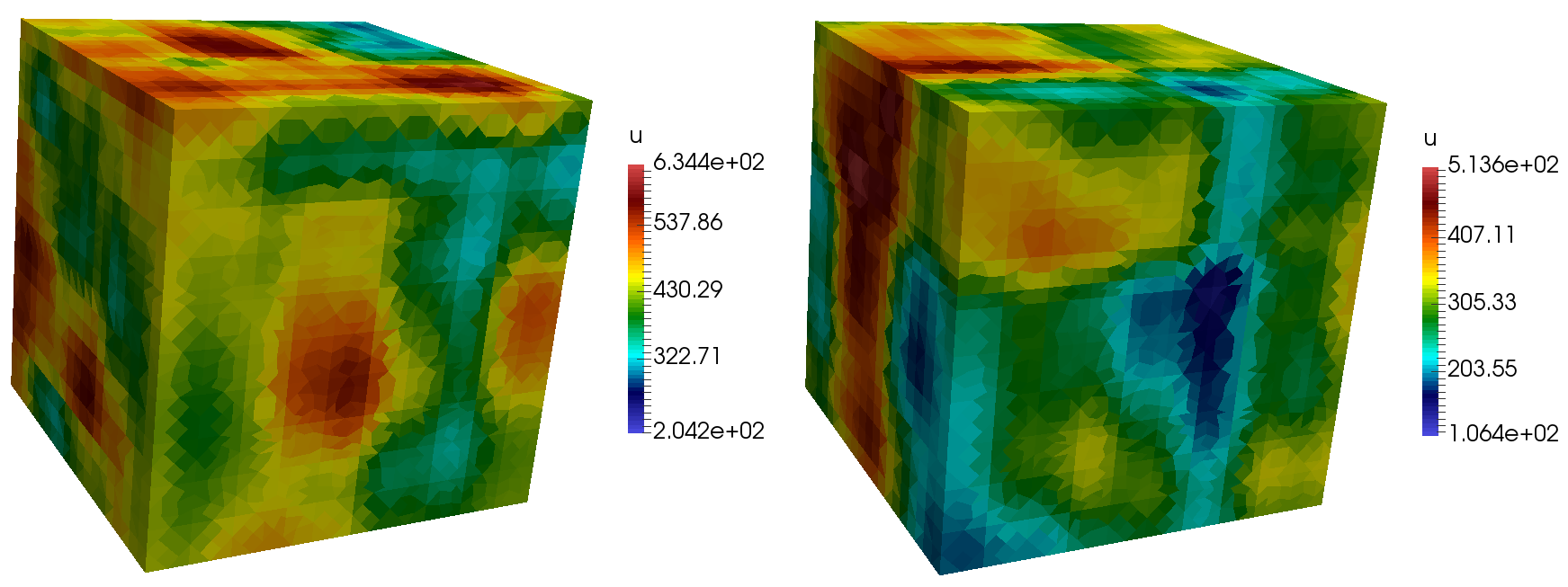} 
\caption{Computational grids and heterogeneous properties for \textit{Test 4} and or \textit{Test 5} (three-dimensional problem). 
First row: coarse grid (blue color),  fine grid (green color) and fracture geometry (red color). 
Second row: Heterogeneous coefficient $k_{1,s}(x)$  (left) and $k_{2,s}(x)$ for \textit{Test 5}  ($k_{m,s}(x)$ for \textit{Test 4}) (right)}
\label{geom3d}
\end{figure}

We present numerical results for several model problems: 
\begin{itemize}
\item[] \textit{Test 1.} Two-dimensional problem with $\Omega = [0,1]^2$,  $p = (p_m, p_f)$ with  $T_{max} = 66 \cdot 10^{-4}$, $k_{m,s} = 10$ and $k_{f,s} = 10^9$ (homogeneous porous matrix).
\item[] \textit{Test 2.} Two-dimensional problem with $\Omega = [0,1]^2$,  $p = (p_m, p_f)$ with $T_{max} = 11 \cdot 10^{-4}$, $k_{m,s} = k_{m,s}(x)$ and $k_{f,s} = 10^9$ (heterogeneous porous matrix).
\item[] \textit{Test 3.} Two-dimensional problem with $\Omega_1 = \Omega_2 = [0,1]^2$,  $p = (p_1, p_2, p_f)$ with $T_{max} = 10^{-4}$, $k_{1,s} = k_{1,s}(x)$, $k_{2,s} = k_{2,s}(x)$ and $k_{f,s} = 10^9$.
\item[] \textit{Test 4.} Three-dimensional problem with $\Omega = [0,10]^3$,  $p = (p_m, p_f)$ with $T_{max}=31 \cdot 10^{-4}$,  $k_{m,s} = 10$ and $k_{f,s} = 10^9$ (homogeneous porous matrix).
\item[] \textit{Test 5.} Three-dimensional problem with $\Omega = [0,10]^3$,  $p = (p_m, p_f)$ with $T_{max}= 1.25 \cdot 10^{-4}$, $k_{m,s} = k_{m,s}(x)$ and $k_{f,s} = 10^9$ (heterogeneous porous matrix).
\item[] \textit{Test 6.} Three-dimensional problem with $\Omega _1 = \Omega_2= [0,10]^3$,  $p = (p_1, p_2, p_f)$ with $T_{max} = 3 \cdot 10^{-4}$, $k_{1,s} = k_{1,s}(x)$, $k_{2,s} = k_{2,s}(x)$ and $k_{f,s} = 10^9$.
\end{itemize}
In Figures \ref{geom2d} for  and \ref{geom3d}, we show computational domains and heterogeneous coefficients for test problems are presented, by green and blue color, we depict the fine and coarse grids. The coarse grids are uniform with $121$ nodes ($100$ cells) for two-dimensional test problems and $726$ nodes ($500$ cells) for three-dimensional tests. 
We use $DOF_f$ (Degree of Freedom) to denote fine grid system size and $DOF_c$ to denote problem size of the multiscale system using GMsFEM.

In Haverkamp model, we set $A_{\alpha} = 1.511 \cdot 10^6$, $B_{\alpha} = 3.96$, $\Theta_{\alpha,s} = 0.287$, $\Theta_{\alpha,r} = 0.075$, $C_{\alpha} = 1.175 \cdot 10^6$ and $D_{\alpha} = 4.74$ for $\alpha = m,f$ for \textit{Tests 1,2,4,5} and $\alpha = 1,2,f$ for \textit{Tests 3, 6}. 
We simulate for $T_{max} = 66 \cdot 10^{-4}$ with 200 time steps and set initial condition $p^0 = -61.5$.  
As boundary conditions, we set $p_g = -20.7$ on the top boundary  ($\Gamma_D$) and zero flux on other boundaries ($\Gamma_N$). 

To compare the results, we use the relative $L^2$ error, relative energy error and relative error for multicontinuum system  between multiscale and fine-scale solutions for continuum $\alpha$  
\[
e^{L_2}_{\alpha} 
=\frac{ ||p_{\alpha} - p_{ms, \alpha} ||_{L_2} }{ ||p_{\alpha} ||_{L_2}}, \quad
e^{a}_{\alpha} 
=\frac{ ||p_{\alpha} - p_{ms, \alpha} ||_{a} }{ ||p_{\alpha}||_{a}}, \quad 
e^{q} 
=\frac{ ||p - p_{ms} ||_{q}}{||p||_{q}}, 
\]
where
\[
|| v_{\alpha} ||_{L_2}^2 = (v_{\alpha}, \, v_{\alpha}), 
\quad 
|| v_{\alpha} ||_{a}^2 = a_{\alpha}(v_{\alpha}, \, v_{\alpha}), 
\quad 
|| v ||_{q}^2 = 
\sum_{\alpha} \left(
a_{\alpha}(v_{\alpha}, \, v_{\alpha}) + 
\sum_{\beta \neq \alpha}
q_{\alpha}(v_{\alpha} - v_{\beta}, \, v_{\alpha})
\right),
\]
and $p_{ms}$ and $p$ are the multiscale and reference (fine-scale) solutions. 
We use GMSH software \cite{gmsh} to construct computational domains and grids. The implementation is based on the open-source library FEniCS \cite{fenics}.

\subsection{Two-dimensional model problems}

We consider the unsaturated fluid flow in two-dimensional computational domain $\Omega = [0,1]^2$  for \textit{Tests 1,2} and \textit{3} (Figure \ref{geom2d}). 
Unstructured fine grid contains approximately 31 409 triangle cells and 14 376 vertices. 
We use $10\times10$ structured coarse grid with rectangular cells to find multiscale solution with $100$ quadratic cells and $121$ vertices .

\begin{figure}[h!]
\centering
\includegraphics[width=0.99\linewidth]{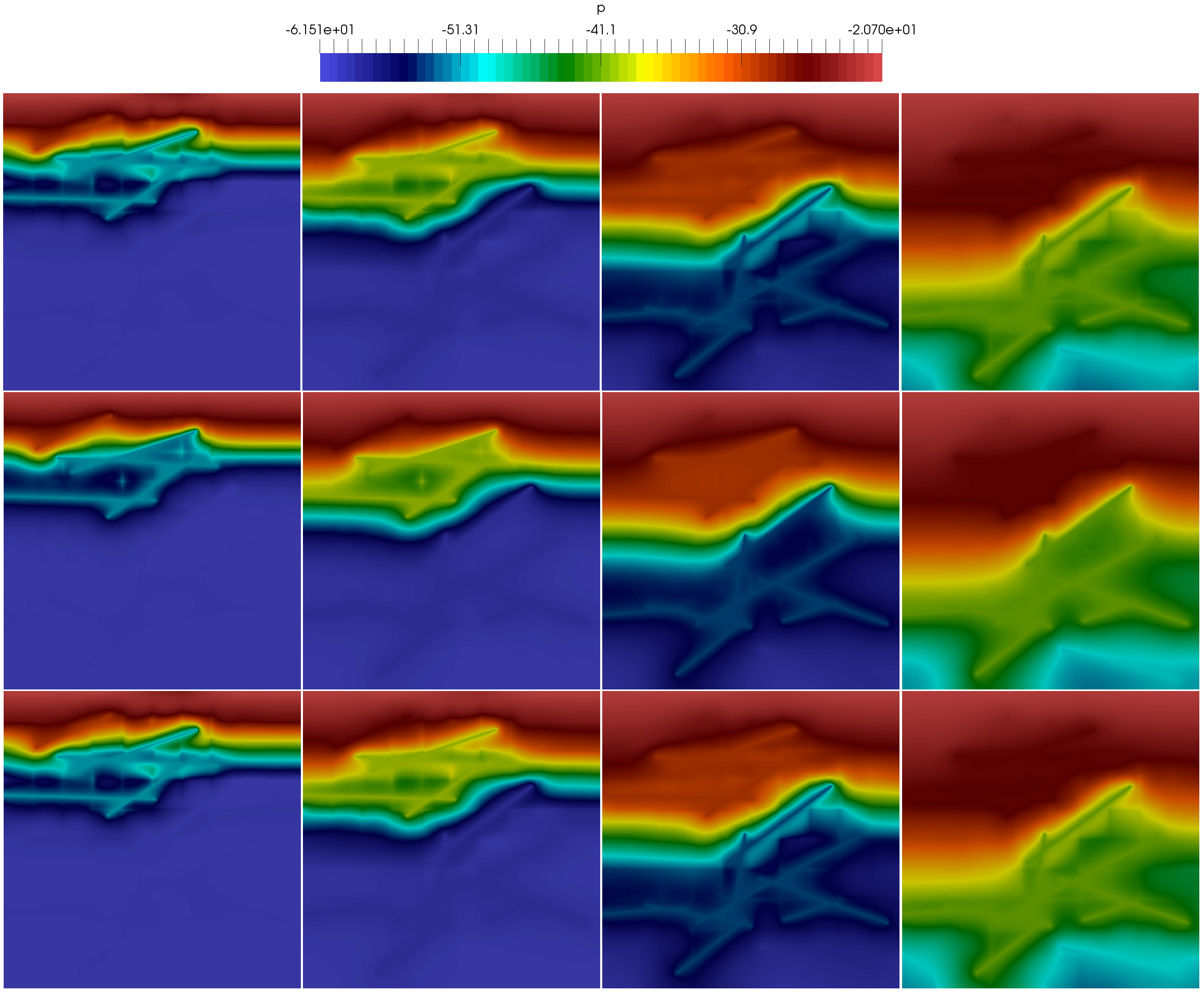} 
\caption{Numerical results for \textit{Test 1}. Solution $p_m$ and $p_{ms,m}$ for different times,  $t_{25}$, $t_{50}$, $t_{100}$ and $t_{200}$ (from left to right). 
First row: fine scale solution $DOF_f=14376$.
Second row: multiscale solution using 8 basis functions $DOF_c=968$. 
Third row: multiscale solution using simplified basis functions $DOF_c=295$}
\label{t1-u}
\end{figure} 

\begin{table}[h!]
\begin{center}
\begin{tabular}{|c|c|cc|cc|cc|cc|}
\hline
Number of multiscale
 &  \multirow{2}{*}{$DOF_c$ }
 & \multicolumn{2}{|c|}{$t_{25}$}& \multicolumn{2}{c|}{$t_{50}$}
 & \multicolumn{2}{|c|}{$t_{100}$} & \multicolumn{2}{c|}{$t_{200}$} \\
basis functions &  
 & $e^{L_2}_m$ & $e^{a}_m$ 
 & $e^{L_2}_m$ & $e^{a}_m$ 
 & $e^{L_2}_m$ & $e^{a}_m$ 
 & $e^{L_2}_m$ & $e^{a}_m$  \\ 
\hline
1 & 121 & 87.47 & 1605.09 & 87.44 & 2243.09 & 87.67 & 2064.03 & 87.87 & 2705.61 \\
2 & 242 & 48.06 & 1836.22 & 43.85 & 2406.58 & 40.39 & 2140.27 & 33.01 & 1989.89 \\
4 & 484 & 1.81 & 24.41 & 1.09 & 12.41 & 0.99 & 9.16 & 0.84 & 8.61 \\
8 & 968 & 0.31 & 3.16 & 0.21 & 1.85 & 0.18 & 1.25 & 0.15 & 1.21 \\
12 & 1452 & 0.15 & 1.35 & 0.09 & 0.69 & 0.09 & 0.62 & 0.08 & 0.56 \\
16 & 1936 & 0.12 & 1.05 & 0.08 & 0.54 & 0.05 & 0.26 & 0.03 & 0.16 \\
\hline
 Simplified & 295 & 1.73 & 28.23 & 1.06 & 17.72 & 1.19 & 17.67 & 1.01 & 17.11 \\
 Adaptive by $\lambda$ & 290 & 2.48 & 30.95 & 1.54 & 18.74 & 2.13 & 24.47 & 1.81 & 24.35 \\ \hline
\end{tabular}
\end{center}
\caption{Numerical results for \textit{Test 1}. Relative $L_2$ and energy errors ($\%$) for different number of multiscale basis functions}
\label{t1-err}
\end{table}

In \textit{Test 1}, we consider problem for $p = (p_m, p_f)$ with homogeneous porous matrix and fracture properties.
The fine-scale solution and multiscale solution using 8 basis function are presented on Figure \ref{t1-u} for four time steps. Relative $L_2$ and energy errors are presented for different number of multiscale basis functions in Table \ref{t1-err}.
We presented error comparison for fixed number of multiscale basis functions in each local domain $\omega_i$ ($M_i = M$ for $i = 1,...,N^{\Omega}_c$) and investigate the adaptive approach that based on the choosing of the number of basis functions that related to the smallest eigenvalues. This adaptive approach, automatically identify number of basis functions that sufficient for obtaining a good results. Furthermore, we show the results for simplified multiscale basis functions that can be applied for some simple cases when we a priory can define main local modes of the solution. We observe that the 1 and 2 basis functions are not sufficient for obtaining a good results, but 4 multiscale basis functions can be sufficient for given geometry.
Presented results demonstrates the dynamic of the multiscale solution behavior with different number of basis functions.
We observe, a good convergence of the solution, when we took sufficient number of multiscale basis functions and with increasing number of basis functions.

\begin{figure}[h!]
\centering
\includegraphics[width=0.99\linewidth]{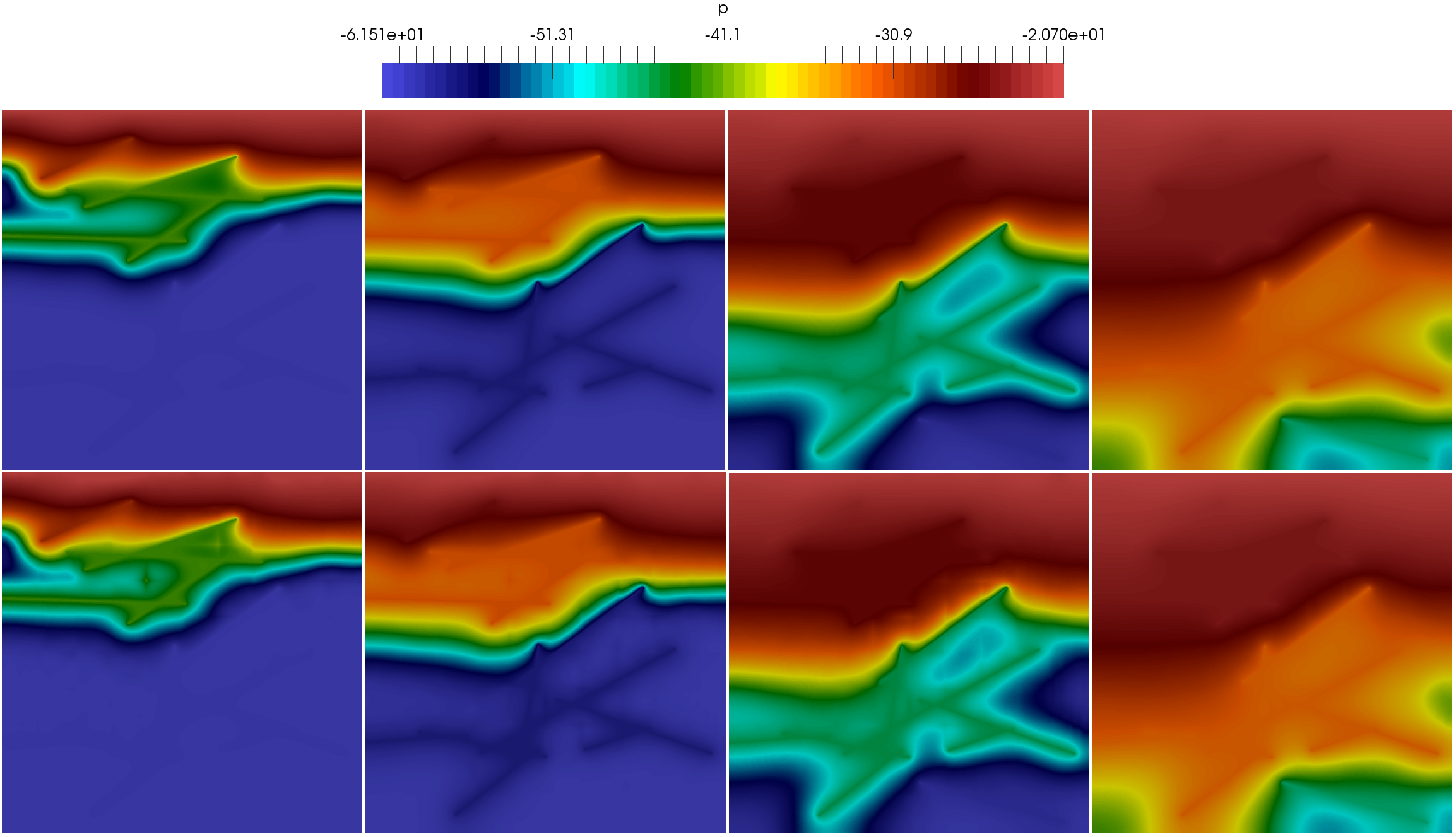} 
\caption{Numerical results for \textit{Test 2}. Solution $p_m$ and $p_{ms,m}$ for different times,  $t_{25}$, $t_{50}$, $t_{100}$ and $t_{200}$ (from left to right). 
First row: fine scale solution $DOF_f=14376$.
Second row: multiscale solution using 8 basis functions $DOF_c=968$}
\label{t2-u}
\end{figure} 

\begin{table}[h!]
\begin{center}
\begin{tabular}{|c|c|cc|cc|cc|cc|}
\hline
Number of multiscale
 &  \multirow{2}{*}{$DOF_c$ }
 & \multicolumn{2}{|c|}{$t_{25}$}& \multicolumn{2}{c|}{$t_{50}$}
 & \multicolumn{2}{|c|}{$t_{100}$} & \multicolumn{2}{c|}{$t_{200}$} \\
basis functions &  
 & $e^{L_2}_m$ & $e^{a}_m$ 
 & $e^{L_2}_m$ & $e^{a}_m$ 
 & $e^{L_2}_m$ & $e^{a}_m$ 
 & $e^{L_2}_m$ & $e^{a}_m$  \\ 
\hline
1 & 121 & 77.37 & 1204.73 & 76.57 & 1418.63 & 74.58 & 1087.39 & 67.86 & 1726.95 \\
2 & 242 & 13.27 & 187.27 & 19.76 & 450.39 & 17.92 & 318.65 & 10.41 & 212.56 \\
4 & 484 & 1.96 & 46.45 & 1.05 & 32.52 & 1.44 & 38.32 & 0.57 & 23.98 \\
8 & 968 & 0.37 & 3.54 & 0.23 & 1.73 & 0.36 & 3.23 & 0.11 & 0.56 \\
12 & 1452 & 0.23 & 1.31 & 0.14 & 0.88 & 0.21 & 1.81 & 0.05 & 0.24 \\
16 & 1936 & 0.13 & 0.91 & 0.07 & 0.39 & 0.11 & 0.78 & 0.02 & 0.09 \\
\hline
\end{tabular}
\end{center}
\caption{Numerical results for \textit{Test 2}. Relative $L_2$ and energy errors ($\%$) for different number of multiscale basis functions}
\label{t2-err}
\end{table}

Next, we consider \textit{Test 2}  for $p = (p_m, p_f)$ with  heterogeneous matrix property $k_{m,s}(x)$ (see Figure \ref{geom2d}). 
The fine-scale solution and multiscale solution using 8 basis function are  presented on Figure \ref{t2-u}. In Table \ref{t2-err}, we present a dynamic of the $L_2$ and energy errors for different number of multiscale basis functions for four time layers.  We observe similarly good results as for previous homogeneous porous matrix property. The GMsFEM demonstrate a good convergence with at least 4 multiscale basis functions.

\begin{figure}[h!]
\centering
\includegraphics[width=0.98\linewidth]{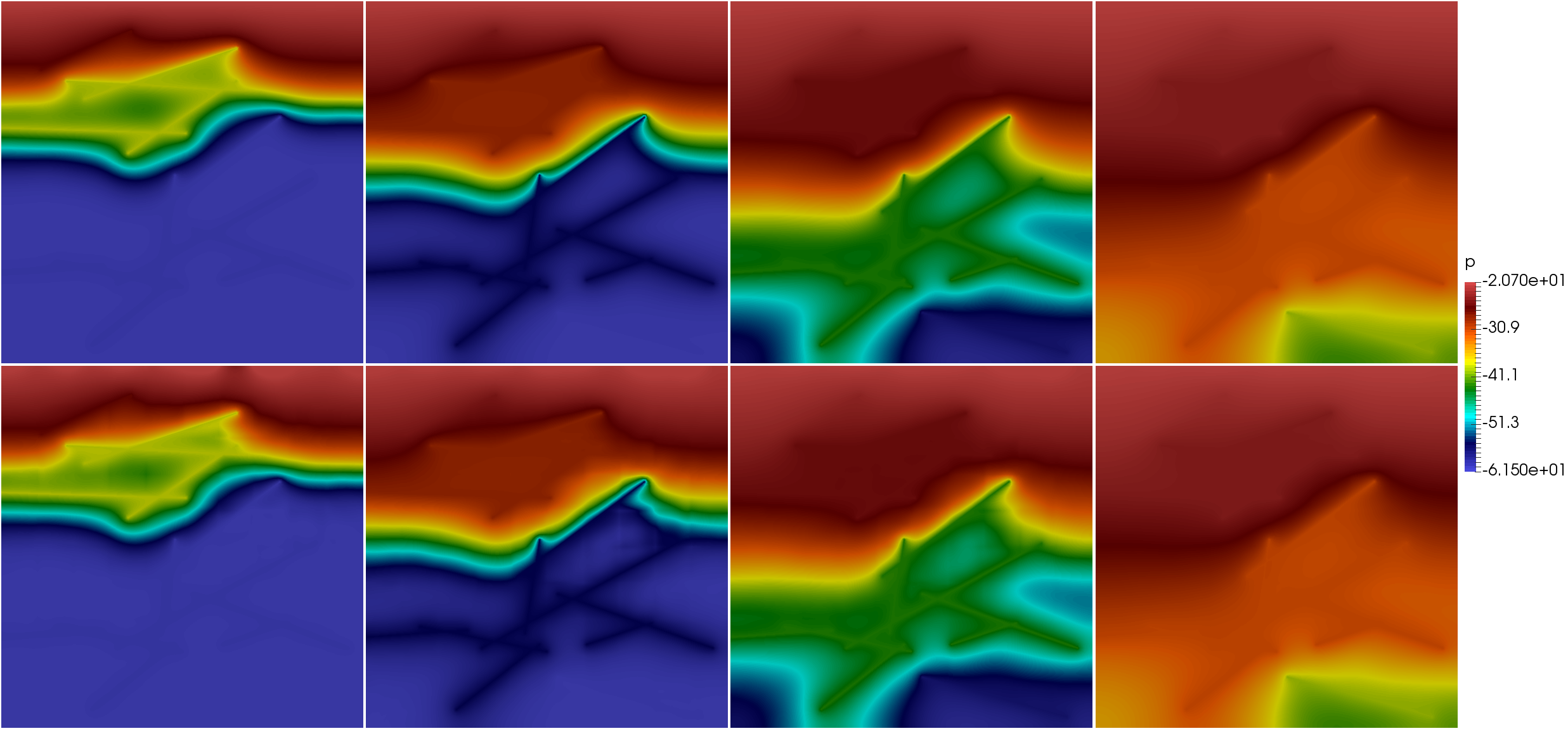} \\
\includegraphics[width=0.98\linewidth]{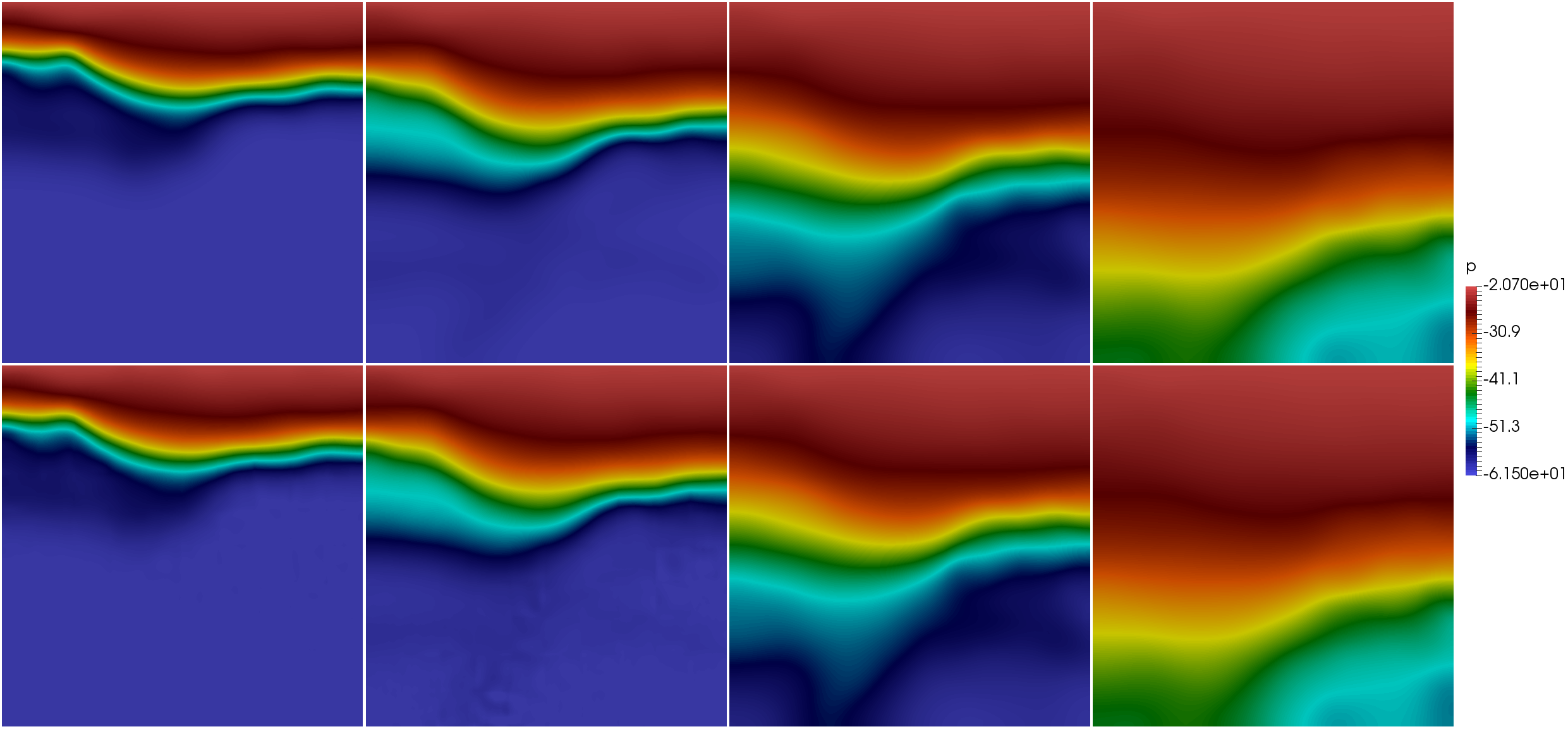} 
\caption{Numerical results for \textit{Test 3}. 
Solution $p_1, p_2$ and $p_{ms,1}, p_{ms,2}$ for different times,  $t_{25}$, $t_{50}$, $t_{100}$ and $t_{200}$ (from left to right). 
First row: fine scale solution. (First continuum).
Second row: multiscale solution using 8 basis functions (First continuum).
Third row: fine scale solution (Second continuum).
Fourth row: multiscale solution using 8 basis functions (Second continuum). $DOF_f=28752$, $DOF_c=968$}
\label{t3-u}
\end{figure} 


\begin{table}[h!]
\begin{center}
\begin{tabular}{|c|c|cc|cc|c|}
\hline 
Number of multiscale
 &  \multirow{2}{*}{$DOF_c$ }
 & \multicolumn{2}{|c|}{$p_1$} & \multicolumn{2}{c|}{$p_2$} & p \\
basis functions &  
 & $e^{L_2}_1$ & $e^{a}_1$  
 & $e^{L_2}_2$ & $e^{a}_2$  & $e^{q}$  \\ 
\hline
1 & 242 & 47.51 & 1228.31 & 0.46 & 2.71 & 1301.74 \\
2 & 484 & 9.12 & 385.99 & 0.27 & 1.98 & 63.74 \\
4 & 968 & 0.71 & 6.47 & 0.11 & 0.62 & 0.39 \\
8 & 1936 & 0.09 & 0.57 & 0.02 & 0.04 & 0.02 \\
12 & 2904 & 0.03 & 0.11 & 0.01 & 0.01 & 0.01 \\
16 & 3872 & 0.13 & 2.85 & 0.01 & 0.02 & 0.04 \\
\hline
\end{tabular}
\end{center}
\caption{Numerical results for \textit{Test 3}. Relative errors ($\%$) for different number of multiscale basis functions}
\label{t3-err}
\end{table}

Next, we consider \textit{Test 3} for $p = (p_1, p_2, p_f)$ with heterogeneous properties $k_{1,s}(x)$ and $k_{2,s}(x)$ with $\sigma _{12} = \sigma _{21} = k_{2,s}(x)$. Reference (fine-scale) and multiscale solutions are presented in Figure \ref{t3-u}. In Table \ref{t3-err}, we present a relative errors and observe good convergence of the presented method with coupled multiscale basis functions for multicontinua unsaturated flow problem.
For example, when we take 4 multiscale basis functions we have less than one percent of $L_2$ errors for $p_1$ and $p_2$ with near 5 \% of energy error.

For two-dimensional test problems (\textit{Test 1, 2} and \textit{3}), we observe small errors, when we take at least 4 multiscale basis functions with huge reduction of the system size. For \textit{Test 1} and \textit{2}, the fine grid system has size $DOF_f = 14 376$, but using coarse grid multiscale solver with 4 basis functions, we reduce size to $DOF_c = 484$. For \textit{Test 3}, fine grid and coarse grid systems have size $DOF_f = 28 752$ and $DOF_C = 968$ for 4 multiscale basis functions, respectively.

\subsection{Three-dimensional model problems}

In this section, we consider three-dimensional model problems and present a numerical results for \textit{Test 4,5} and \textit{6}. We consider, similar test problems as for two - dimensional problems.
We start with \textit{Test 4} with homogeneous porous matrix properties and fracture distribution that presented on Figure \ref{geom3d}.
Heterogeneous coefficient for \textit{Test 4} and \textit{5} are presented in Figure \ref{geom3d}. Fine grid has $125 740$ tetrahedral cells and $21 609$ vertices. In GMsFEM, we use $5 \times 5 \times 5$ structured coarse grid with $125$ cubic cells and $216$ vertices.

\begin{figure}[h!]
\centering
\includegraphics[width=0.8\linewidth]{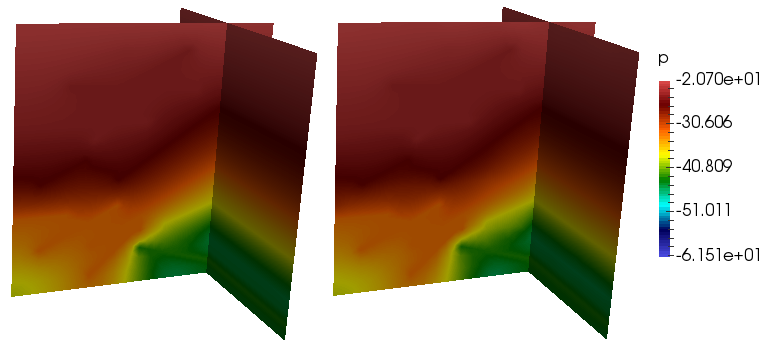} 
\caption{Numerical results for \textit{Test 4}. 
Solution $p_m$ and $p_{ms,m}$ for final time.  
Left: fine scale solution $DOF_f=21609$.
Right: multiscale solution using 16 basis functions $DOF_c=3456$}
\label{t4-u}
\end{figure}

\begin{table}[h!]
\begin{center}
\begin{tabular}{|c|c|cc|cc|cc|cc|}
\hline
Number of multiscale
 &  \multirow{2}{*}{$DOF_c$ }
 & \multicolumn{2}{|c|}{$t_{25}$}& \multicolumn{2}{c|}{$t_{50}$}
 & \multicolumn{2}{|c|}{$t_{100}$} & \multicolumn{2}{c|}{$t_{200}$} \\
basis functions &  
 & $e^{L_2}_m$ & $e^{a}_m$ 
 & $e^{L_2}_m$ & $e^{a}_m$ 
 & $e^{L_2}_m$ & $e^{a}_m$ 
 & $e^{L_2}_m$ & $e^{a}_m$  \\ 
\hline
1 & 216 & 29.68 & 102.99 & 36.41 & 101.46 & 37.23 & 103.83 & 17.39 & 108.95 \\
2 & 432 & 22.04 & 194.11 & 31.17 & 286.69 & 31.23 & 292.51 & 13.91 & 240.97 \\
4 & 864 & 6.73 & 89.81 & 2.21 & 13.91 & 4.59 & 44.39 & 0.95 & 10.79 \\
8 & 1728 & 5.07 & 61.19 & 2.21 & 13.33 & 2.09 & 18.59 & 0.74 & 7.05 \\
12 & 2592 & 2.32 & 18.16 & 0.73 & 4.74 & 1.05 & 7.33 & 0.31 & 2.78 \\
16 & 3456 & 1.76 & 14.35 & 0.68 & 4.26 & 0.73 & 4.19 & 0.21 & 1.77 \\
\hline
\end{tabular}
\end{center}
\caption{Numerical results for \textit{Test 4}. Relative errors ($\%$) for different number of multiscale basis functions}
\label{t4-err}
\end{table}

\begin{figure}[h!]
\centering
\includegraphics[width=0.8\linewidth]{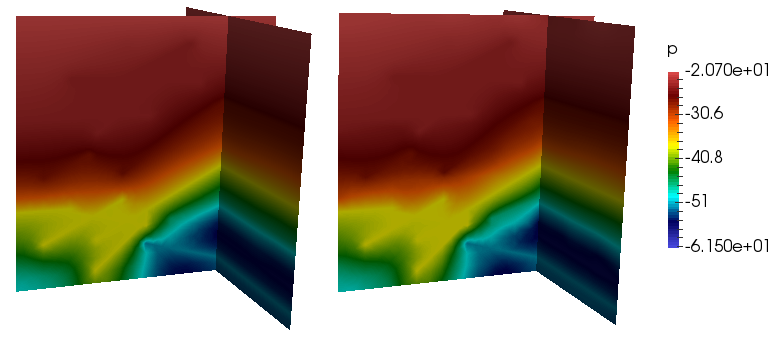} 
\caption{Numerical results for \textit{Test 5}. 
Solution $p_m$ and $p_{ms,m}$ for final time.  
Left: fine scale solution $DOF_f=21609$.
Right: multiscale solution using 16 basis functions $DOF_c=3456$}
\label{t5-u}
\end{figure} 

\begin{table}[h!]
\begin{center}
\begin{tabular}{|c|c|cc|cc|cc|cc|}
\hline
Number of multiscale
 &  \multirow{2}{*}{$DOF_c$ }
 & \multicolumn{2}{|c|}{$t_{25}$}& \multicolumn{2}{c|}{$t_{50}$}
 & \multicolumn{2}{|c|}{$t_{100}$} & \multicolumn{2}{c|}{$t_{200}$} \\
basis functions &  
 & $e^{L_2}_m$ & $e^{a}_m$ 
 & $e^{L_2}_m$ & $e^{a}_m$ 
 & $e^{L_2}_m$ & $e^{a}_m$ 
 & $e^{L_2}_m$ & $e^{a}_m$  \\ 
\hline
1 & 216 & 27.64 & 105.15 & 34.61 & 100.73 & 38.99 & 103.66 & 28.34 & 106.47 \\
2 & 432 & 18.97 & 174.97 & 28.97 & 260.03 & 33.77 & 311.93 & 22.62 & 264.14 \\
4 & 864 & 6.82 & 86.51 & 2.72 & 23.15 & 3.37 & 31.11 & 2.38 & 21.24 \\
8 & 1728 & 7.61 & 116.12 & 3.17 & 28.67 & 2.16 & 18.87 & 1.23 & 10.15 \\
12 & 2592 & 3.97 & 39.72 & 1.36 & 10.52 & 1.13 & 8.79 & 0.72 & 5.61 \\
16 & 3456 & 4.12 & 43.41 & 1.46 & 10.41 & 0.92 & 5.99 & 0.55 & 3.83 \\
\hline
\end{tabular}
\end{center}
\caption{Numerical results for \textit{Test 5}. Relative errors ($\%$) for different number of multiscale basis functions}
\label{t5-err}
\end{table}

The fine-scale solution and multiscale solution using 16 basis function are presented on Figures \ref{t4-u} and \ref{t5-u} for \textit{Test 4} and \textit{5}, respectively. In Tables \ref{t4-err} and \ref{t5-err}, we present a dynamic for the relative $L_2$ and energy errors for \textit{Test 4} and \textit{5}. From the results we can conclude that this method has a good performance for three-dimensional problems with heterogeneous and homogeneous porous matrix properties. For example, when we take 16 multiscale basis functions, we obtain solution with $0.21 \%$ of relative $L_2$ error and $1.77 \%$ of energy error at final time for \textit{Test 4}. Fine grid system has size $DOF_f = 21 609$ and $DOF_c = 3456$ for multiscale solver with 16 basis functions.

\begin{figure}[h!]
\centering
\includegraphics[width=0.8\linewidth]{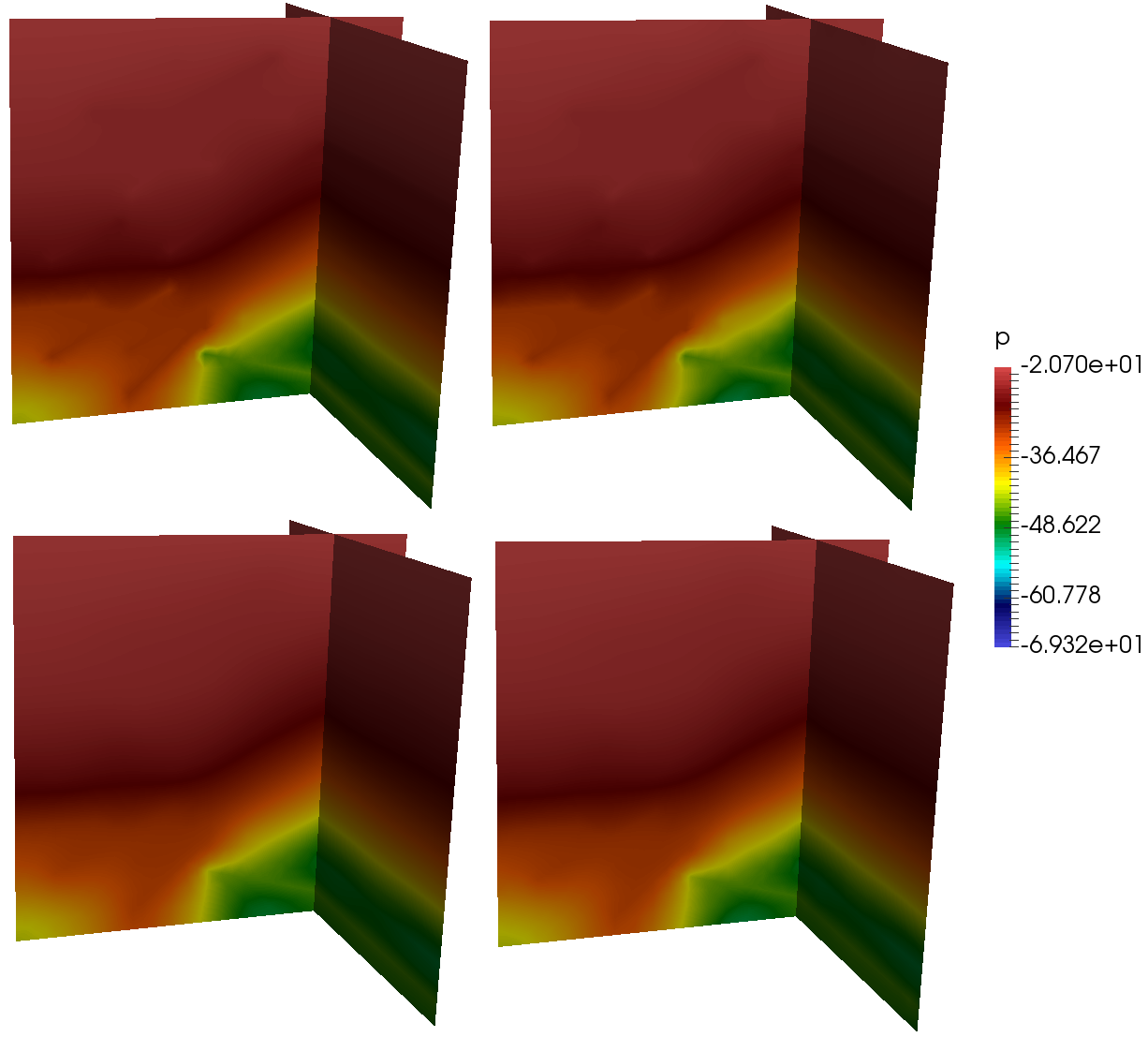} 
\caption{Numerical results for \textit{Test 6}. 
Solution $p_1, p_2$ and $p_{ms,1}, p_{ms,2}$ for final time.  
Left: fine scale solution $DOF_f=43218$.
Right: multiscale solution using 16 basis functions $DOF_c=6912$}
\label{t6-u}
\end{figure} 

\begin{table}[h!]
\begin{center}
\begin{tabular}{|c|c|cc|cc|c|}
\hline 
Number of multiscale
 &  \multirow{2}{*}{$DOF_c$ }
 & \multicolumn{2}{|c|}{$p_1$} & \multicolumn{2}{c|}{$p_2$} & p \\
basis functions &  
 & $e^{L_2}_1$ & $e^{a}_1$  
 & $e^{L_2}_2$ & $e^{a}_2$  & $e^{q}$  \\ 
\hline
2 & 864 & 21.77 & 263.19 & 12.88 & 111.66 & 2849.94 \\
4 & 1728 & 2.23 & 19.44 & 1.17 & 7.56 & 3.64 \\
8 & 3456 & 1.18 & 9.89 & 0.68 & 3.74 & 0.29 \\
12 & 5184 & 0.65 & 5.53 & 0.34 & 1.91 & 0.12 \\
16 & 6912 & 0.49 & 3.73 & 0.24 & 1.15 & 0.08 \\
\hline
\end{tabular}
\end{center}
\caption{Numerical results for \textit{Test 6}. Relative errors ($\%$) for different number of multiscale basis functions}
\label{t6-err}
\end{table}

Fine grid and multiscale solutions for \textit{Test 6} with $p = (p_1, p_2, p_f)$ are presented in Figure \ref{t6-u} for final time. Relative errors are presented in Table \ref{t6-err} for final time. Coefficient $\sigma _{12} = \sigma _{21} = 10^3 \cdot k_{2,s}(x)$. Fine grid system has size $DOF_f = 43218$ and $DOF_c = 6912$ for multiscale solver with 16 basis functions.
We observe a good convergence of the presented method with coupled multiscale basis functions for unsaturated multicontinua flow problems in three - dimensional formulations.

\section{Conclusion}

We presented a multiscale method for simulations of the multicontinua unsaturated flow problems in heterogeneous fractured porous media.
The mathematical model is presented as general multicontinua model. The fine grid approximation are presented for coupled system of equations based on the finite element approximation and Discrete Fracture Model with unstructured grids that resolve fracture geometries explicitly.
We presented the multiscale solver using GMsFEM with coupled multiscale basis function construction for multicontinua problems. Adaptive approach is investigated with simplified multiscale basis functions.
To illustrate the idea of our approach, we considered unsaturated flow in fractured porous media and dual continua background model with discrete fractures networks. Numerical results are presented for two and three dimensional test problems. Presented multiscale solver accurately capture detailed interactions between multiple continua and provide a good results with huge reduction of the discrete system size.
The comparison of the relative error for different number of basis functions and for adaptive approaches are presented for six test problems.

\section{Acknowledgements}
MV's and DS's works are supported by the mega-grant of the Russian Federation Government N14.Y26.31.0013 and RSF N17-71-20055.
The research of Eric Chung is partially supported by the Hong Kong RGC General Research Fund (Project numbers 14304217 and 14302018) and CUHK Faculty of Science Direct Grant 2017-18.

\bibliographystyle{unsrt}
\bibliography{lit}

\begin{thebibliography}{10}

\bibitem{b1}
Michael~A Celia, Efthimios~T Bouloutas, and Rebecca~L Zarba.
\newblock A general mass-conservative numerical solution for the unsaturated
  flow equation.
\newblock {\em Water resources research}, 26(7):1483--1496, 1990.

\bibitem{b2}
Michael~A Celia and Philip Binning.
\newblock A mass conservative numerical solution for two-phase flow in porous
  media with application to unsaturated flow.
\newblock {\em Water Resources Research}, 28(10):2819--2828, 1992.

\bibitem{b4}
Roland Haverkamp, Michel Vauclin, Jaoudat Touma, PJ~Wierenga, and Georges
  Vachaud.
\newblock A comparison of numerical simulation models for one-dimensional
  infiltration 1.
\newblock {\em Soil Science Society of America Journal}, 41(2):285--294, 1977.

\bibitem{lee2001hierarchical}
Seong~H Lee, MF~Lough, and CL~Jensen.
\newblock Hierarchical modeling of flow in naturally fractured formations with
  multiple length scales.
\newblock {\em Water resources research}, 37(3):443--455, 2001.

\bibitem{li2006efficient}
Liyong Li, Seong~Hee Lee, et~al.
\newblock Efficient field-scale simulation for black oil in a naturally
  fractured reservoir via discrete fracture networks and homogenized media.
\newblock In {\em International oil \& gas conference and exhibition in China}.
  Society of Petroleum Engineers, 2006.

\bibitem{chung2017coupling}
Eric~T Chung, Yalchin Efendiev, Tat Leung, and Maria Vasilyeva.
\newblock Coupling of multiscale and multi-continuum approaches.
\newblock {\em GEM-International Journal on Geomathematics}, 8(1):9--41, 2017.

\bibitem{li2018multiscale}
Qiuqi Li, Yuhe Wang, and Maria Vasilyeva.
\newblock Multiscale model reduction for fluid infiltration simulation through
  dual-continuum porous media with localized uncertainties.
\newblock {\em Journal of Computational and Applied Mathematics}, 336:127--146,
  2018.

\bibitem{b5}
Yalchin Efendiev and Thomas~Y Hou.
\newblock Multiscale finite element methods: theory and applications.
\newblock 4, 2009.

\bibitem{hajibeygi2008iterative}
Hadi Hajibeygi, Giuseppe Bonfigli, Marc~Andre Hesse, and Patrick Jenny.
\newblock Iterative multiscale finite-volume method.
\newblock {\em Journal of Computational Physics}, 227(19):8604--8621, 2008.

\bibitem{weinan2007heterogeneous}
E~Weinan, Bjorn Engquist, Xiantao Li, Weiqing Ren, and Eric Vanden-Eijnden.
\newblock Heterogeneous multiscale methods: a review.
\newblock {\em Commun. Comput. Phys}, 2(3):367--450, 2007.

\bibitem{CYH2016adaptive}
Eric Chung, Yalchin Efendiev, and Thomas~Y Hou.
\newblock Adaptive multiscale model reduction with generalized multiscale
  finite element methods.
\newblock {\em Journal of Computational Physics}, 320:69--95, 2016.

\bibitem{chung2018constraint}
Eric~T Chung, Yalchin Efendiev, and Wing~Tat Leung.
\newblock Constraint energy minimizing generalized multiscale finite element
  method.
\newblock {\em Computer Methods in Applied Mechanics and Engineering},
  339:298--319, 2018.

\bibitem{he2009adaptive}
Xinguang He and Li~Ren.
\newblock An adaptive multiscale finite element method for unsaturated flow
  problems in heterogeneous porous media.
\newblock {\em Journal of hydrology}, 374(1-2):56--70, 2009.

\bibitem{he2006multiscale}
Xinguang He and Li~Ren.
\newblock A multiscale finite element linearization scheme for the unsaturated
  flow problems in heterogeneous porous media.
\newblock {\em Water resources research}, 42(8), 2006.

\bibitem{b3}
Victor~Eralingga Ginting.
\newblock Computational upscaled modeling of heterogeneous porous media flow
  utilizing finite volume method.
\newblock 2005.

\bibitem{chen2005upscaling}
Zhiming Chen, Weibing Deng, and Huang Ye.
\newblock Upscaling of a class of nonlinear parabolic equations for the flow
  transport in heterogeneous porous media.
\newblock {\em Communications in Mathematical Sciences}, 3(4):493--515, 2005.

\bibitem{akkutlu2015multiscale}
IY~Akkutlu, Yalchin Efendiev, and Maria Vasilyeva.
\newblock Multiscale model reduction for shale gas transport in fractured
  media.
\newblock {\em Computational Geosciences}, pages 1--21, 2015.

\bibitem{ctene2016algebraic}
Matei {\c{T}}ene, Mohammed~Saad Al~Kobaisi, and Hadi Hajibeygi.
\newblock Algebraic multiscale method for flow in heterogeneous porous media
  with embedded discrete fractures (f-ams).
\newblock {\em Journal of Computational Physics}, 321:819--845, 2016.

\bibitem{bosma2017multiscale}
Sebastian Bosma, Hadi Hajibeygi, Matei Tene, and Hamdi~A Tchelepi.
\newblock Multiscale finite volume method for discrete fracture modeling on
  unstructured grids (ms-dfm).
\newblock {\em Journal of Computational Physics}, 2017.

\bibitem{ctene2017projection}
Matei {\c{T}}ene, Sebastian~BM Bosma, Mohammed~Saad Al~Kobaisi, and Hadi
  Hajibeygi.
\newblock Projection-based embedded discrete fracture model (pedfm).
\newblock {\em Advances in Water Resources}, 105:205--216, 2017.

\bibitem{hkj12}
H.~Hajibeygi, D.~Kavounis, and P.~Jenny.
\newblock A hierarchical fracture model for the iterative multiscale finite
  volume method.
\newblock {\em Journal of Computational Physics}, 230(24):8729--8743, 2011.

\bibitem{akkutlu2018multiscale}
I~Yucel Akkutlu, Yalchin Efendiev, Maria Vasilyeva, and Yuhe Wang.
\newblock Multiscale model reduction for shale gas transport in poroelastic
  fractured media.
\newblock {\em Journal of Computational Physics}, 353:356--376, 2018.

\bibitem{chung2018non}
Eric~T Chung, Yalchin Efendiev, Wing~Tat Leung, Maria Vasilyeva, and Yating
  Wang.
\newblock Non-local multi-continua upscaling for flows in heterogeneous
  fractured media.
\newblock {\em Journal of Computational Physics}, 372:22--34, 2018.

\bibitem{vasilyeva2019nonlocal}
Maria Vasilyeva, Eric~T Chung, Siu~Wun Cheung, Yating Wang, and Georgy
  Prokopev.
\newblock Nonlocal multicontinua upscaling for multicontinua flow problems in
  fractured porous media.
\newblock {\em Journal of Computational and Applied Mathematics}, 2019.

\bibitem{vasilyeva2019constrained}
Maria Vasilyeva, Eric~T Chung, Yalchin Efendiev, and Jihoon Kim.
\newblock Constrained energy minimization based upscaling for coupled flow and
  mechanics.
\newblock {\em Journal of Computational Physics}, 376:660--674, 2019.

\bibitem{b6}
I~Yucel Akkutlu, Yalchin Efendiev, Maria Vasilyeva, and Yuhe Wang.
\newblock Multiscale model reduction for shale gas transport in poroelastic
  fractured media.
\newblock {\em Journal of Computational Physics}, 353:356--376, 2018.

\bibitem{b8}
Eric~T Chung, Wing~Tat Leung, Maria Vasilyeva, and Yating Wang.
\newblock Multiscale model reduction for transport and flow problems in
  perforated domains.
\newblock {\em Journal of Computational and Applied Mathematics}, 330:519--535,
  2018.

\bibitem{b10}
Eric~T Chung, Yalchin Efendiev, Guanglian Li, and Maria Vasilyeva.
\newblock Generalized multiscale finite element methods for problems in
  perforated heterogeneous domains.
\newblock {\em Applicable Analysis}, 95(10):2254--2279, 2016.

\bibitem{b11}
Yalchin Efendiev, Thomas~Y Hou, Victor Ginting, et~al.
\newblock Multiscale finite element methods for nonlinear problems and their
  applications.
\newblock {\em Communications in Mathematical Sciences}, 2(4):553--589, 2004.

\bibitem{b12}
Yalchin Efendiev, Juan Galvis, and Thomas~Y Hou.
\newblock Generalized multiscale finite element methods (gmsfem).
\newblock {\em Journal of Computational Physics}, 251:116--135, 2013.

\bibitem{martin2005modeling}
Vincent Martin, J{\'e}r{\^o}me Jaffr{\'e}, and Jean~E Roberts.
\newblock Modeling fractures and barriers as interfaces for flow in porous
  media.
\newblock {\em SIAM Journal on Scientific Computing}, 26(5):1667--1691, 2005.

\bibitem{d2012mixed}
Carlo D’Angelo and Anna Scotti.
\newblock A mixed finite element method for darcy flow in fractured porous
  media with non-matching grids.
\newblock {\em ESAIM: Mathematical Modelling and Numerical Analysis},
  46(2):465--489, 2012.

\bibitem{formaggia2014reduced}
Luca Formaggia, Alessio Fumagalli, Anna Scotti, and Paolo Ruffo.
\newblock A reduced model for darcy’s problem in networks of fractures.
\newblock {\em ESAIM: Mathematical Modelling and Numerical Analysis},
  48(4):1089--1116, 2014.

\bibitem{Quarteroni2008coupling}
Carlo D'angelo and Alfio Quarteroni.
\newblock On the coupling of 1d and 3d diffusion-reaction equations:
  application to tissue perfusion problems.
\newblock {\em Mathematical Models and Methods in Applied Sciences},
  18(08):1481--1504, 2008.

\bibitem{efendiev2015hierarchical}
Yalchin Efendiev, Seong Lee, Guanglian Li, Jun Yao, and Na~Zhang.
\newblock Hierarchical multiscale modeling for flows in fractured media using
  generalized multiscale finite element method.
\newblock {\em GEM-International Journal on Geomathematics}, 6(2):141--162,
  2015.

\bibitem{eh09}
Y.~Efendiev and T.~Hou.
\newblock {Multiscale Finite Element Methods: Theory and Applications}.
\newblock 4, 2009.

\bibitem{akkutlu2017multiscale}
I~Yucel Akkutlu, Yalchin Efendiev, Maria Vasilyeva, and Yuhe Wang.
\newblock Multiscale model reduction for shale gas transport in a coupled
  discrete fracture and dual-continuum porous media.
\newblock {\em Journal of Natural Gas Science and Engineering}, 2017.

\bibitem{gmsh}
Software gmsh. (http://geuz.org/gmsh/).

\bibitem{fenics}
Anders Logg, Kent-Andre Mardal, and Garth Wells.
\newblock Automated solution of differential equations by the finite element
  method: The fenics book.
\newblock 84, 2012.

\end{thebibliography}
 
\end{document}